\documentclass[12pt]{article}
\usepackage{a4wide}
\usepackage{amssymb}
\usepackage{amsfonts}
\usepackage{amsmath} 

\input xy
\xyoption{matrix}\xyoption{arrow}
\def\edge{\ar@{-}}
\def\dedge{\ar@{.}}

\long\def\ignore#1{#1}

\newtheorem{proposition}{Proposition}[section]
\newtheorem{theorem}[proposition]{Theorem}
\newtheorem{lemma}[proposition]{Lemma}

\newtheorem{definition}[proposition]{Definition}
\newtheorem{corollary}[proposition]{Corollary}

\newcommand{\smallbox}{{\vrule height3pt width3pt depth0pt}}

\def\l{\left(}
\def\r{\right)}
\def\sse{\subseteq}
\def\goesto{\longrightarrow} 
\def\ot{\otimes}

\def\qdot{(-q)^{\bullet}}

\newcommand{\qs}[1]{\left[#1\right]}
\newcommand{\Q}{\bigskip\par\noindent}

\def\o{{\cal O}}
\def\oq{{\cal O}_q}
\def\oqminus{{\cal O}_{q^{-1}}}
\def\oqmm{\oq(M_{m})}
\def\oqmn{\oq(M_{mn})}
\def\oqnm{\oq(M_{nm})}
\def\oqnn{\oq(M_{n})}
\def\oquu{\oq(M_{u})}

\def\oqminusnn{\oqminus(M_{n})}
\def\oqslm{\oq(SL_m)}
\def\oqglm{\oq(GL_m)}

\def\oqglu{\oq(GL_u)}
\def\oqv{\oq(V)}

\def\gh{{\cal G}_h}
\def\ghmn{\gh(m,n)}
\def\gq{{\cal G}_q}
\def\gqminus{{\cal G}_{q^{-1}}}
\def\gqminusmn{\gqminus(m,n)}
\def\gqmn{\gq(m,n)}
\def\gqtwofour{\gq(2,4)}
\def\gqntwon{\gq(n,2n)}

\def\ce{{\cal E}}
\def\cg{{\cal G}} 
\def\u{\qs{n+1\ldots 2n}}
\def\Itil{\widetilde{I}}
\def\Jtil{\widetilde{J}}
\def\Ltil{\widetilde{L}}

\def\mc{{\mathbb C}} 
\def\mn{{\mathbb N}}
\def\mz{{\mathbb Z}}

\def\bfdelta{{\bf \delta}}

\def\qmij{[I\mid J]}
\def\oqij{\oq(M_{I,J})}

\newcommand{\oqijlong}[1]{[i_1 \dots i_{#1}\mid j_1 \dots j_{#1}]} 
\def\co{{\rm co}}
\def\dhom{{\rm Dhom}}
\def\dhomrx{\dhom(R,x)}
\def\gkdim{{\rm GKdim}}
\def\id{{\rm id}}
\def\lex{<_{\rm{lex}}}
\def\blex{>_{\rm{lex}}}

\begin{document}

\title{Ring theoretic properties of quantum grassmannians}
\author{A C Kelly\footnote{Some of the results in this paper appear in
the first author's PhD thesis (Edinburgh, 2001). She thanks EPSRC for
financial support.}\,, T H Lenagan and L Rigal\footnote{Part of this work
was done while the third author was visiting the
University of Edinburgh. He thanks the Edinburgh Mathematical
Society for the financial support of this visit.}
}
\date{}
\maketitle

\begin{abstract}The $m\times n$ quantum grassmannian, $\gqmn$, with
$m\leq n$,  is the
subalgebra of the algebra ${\cal O}_q(M_{mn})$ 
of quantum $m\times n$ matrices that is
generated by the maximal $m\times m$ quantum minors. Several properties
of $\gqmn$ are established. In particular, a $k$-basis of $\gqmn$ is
obtained, and it is shown that $\gqmn$ is a noetherian domain of
Gelfand-Kirillov dimension $m(n-m) +1$. The algebra $\gqmn$ is
identified as the subalgebra of coinvariants of a natural left coaction
of $\oqslm$ on ${\cal O}_q(M_{mn})$ 
and it is shown that $\gqmn$ is a maximal order.

\Q
\end{abstract}

\section*{Introduction}

\Q
Fix a base field $k$, a nonzero scalar $q\in k$ and positive integers $m,
n$ with $m\leq n$. The {\em coordinate
ring of quantum $m\times n$ matrices}, $\oqmn$, is the $k$-algebra
generated by $mn$ indeterminates $X_{ij}$, $1 \le i \le
m$ and $1 \le j \le n$, subject to the following relations:
\begin{equation} \label{oqmn}
\begin{array}{rcl}
X_{ij}X_{il} &=& qX_{il}X_{ij}, \\
X_{ij}X_{kj} &=& qX_{kj}X_{ij}, \\
X_{il}X_{kj} &=& X_{kj}X_{il},  \\
X_{ij}X_{kl} - X_{kl}X_{ij} &=& (q-q^{-1})X_{il}X_{kj},
\end{array}
\end{equation}
for $1 \le i < k \le m$ and $1 \le j < l \le n$. It is well-known that
$\oqmn$ can be presented as an iterated skew polynomial algebra over $k$
with the generators added in lexicographic order. As a consequence of
this presentation, it is easy to establish that $\oqmn$ is a noetherian
domain of Gelfand-Kirillov dimension $mn$. 

We will usually write $\oqnn$ for the algebra $\oq(M_{nn})$. In this
algebra the {\em quantum determinant}, $D_q = \det_q$ is defined by 
\[
D_q := \sum_{\sigma \in S_{n}}(-q)^{l(\sigma)}
X_{1,\sigma(1)} \dots X_{n,\sigma(n)};
\]
from \cite[Theorem 4.6.1]{PW},
we know that $D_q$ is in the centre of $\oqnn$. 

Following \cite{GL}, we use the notation $\qmij$ to denote the quantum
determinant of the quantum matrix subalgebra $\oqij$ of ${\cal
O}_q(M_{mn})$ 
generated by the
elements $X_{ij}$ with $i\in I$ and $j\in J$, where $I$ and $J$ are
index sets with $|I| = |J|$. The element $[I|J]$ is the {\em quantum
minor} determined by the index sets $I$ and $J$. If $I = \{i_1, \dots, 
i_s\}$ and $J = \{ j_1, \dots, j_s\}$ where the indices are written in
ascending order, then we will often denote $\qmij$ by $\oqijlong{s}$. 

In this paper we are interested in studying the ring theoretic
properties of a certain subalgebra of $\oqmn$, the {\em 
quantum deformation of the homogeneous 
coordinate ring of the $m\times n$
grassmannian}, $\gqmn$.  This is a
deformation of the classical homogeneous 
coordinate ring of the grassmannian of
$m$-dimensional $k$-subspaces of $n$-dimensional $k$-space and is
generated by the maximal quantum minors of $\oqmn$; to be more specific,
$\gqmn$ is the subalgebra of $\oqmn$ generated by the $m\times m$
quantum minors of $\oqmn$.  In the quantum grassmannian $\gqmn$, any
$m\times m$ quantum minor will involve rows $1, \dots, m$ of the quantum
matrix $(X_{ij})$ associated to 
$\oqmn$.  Thus, to simplify notation, we may denote a quantum
minor by its columns only; that is, the quantum minor given by the row
set $\{1, \dots, m\}$ and column set $J$ will be denoted by $[J]$. 

\Q{\bf Example}
$\gqtwofour$ is the $k$-algebra generated by the $2\times 2$ minors of
the $2
\times 4$ quantum matrix of ${\cal
O}_q(M_{2,4})$: $\qs{12}, \qs{13} , \qs{14}, \qs{23},
\qs{24}$ and  $\qs{34} .$\\
Using the relations for $\oqmn$ and \cite[Lemma A.1]{GL} 
we can  calculate the following commutation
relations:
\[
\qs{12}\qs{13}  =  q\qs{13}\qs{12},\;\;\;
\qs{12}\qs{14}  =  q\qs{14}\qs{12},\;\;\;
\qs{12}\qs{23}  =  q\qs{23}\qs{12},\]\[
\qs{12}\qs{24}  =  q\qs{24}\qs{12},\;\;\;
\qs{12}\qs{34}  =  q^2\qs{34}\qs{12},\;\;\;
\qs{13}\qs{14}  =  q\qs{14}\qs{13},\]\[  
\qs{13}\qs{23}  =  q\qs{23}\qs{13},\;\;\;
\qs{13}\qs{24}  =  \qs{24}\qs{13} + \l q-q^{-1} \r \qs{14}\qs{23},\]\[
\qs{13}\qs{34}  =  q\qs{34}\qs{13},\;\;\;  
\qs{14}\qs{23}  =  \qs{23}\qs{14},\;\;\; 
\qs{14}\qs{24}  =  q\qs{24}\qs{14},\]\[  
\qs{14}\qs{34}  =  q\qs{34}\qs{14},\;\;\;
\qs{23}\qs{24}  =  q\qs{24}\qs{23},\;\;\;  
\qs{23}\qs{34}  =  q\qs{34}\qs{23},\]\[ 
\qs{24}\qs{34}  =  q\qs{34}\qs{24},
\]
and the Quantum Pl\"{u}cker relation
\[
\qs{12}\qs{34} - q\qs{13}\qs{24} +q^2\qs{14}\qs{23} = 0.
\]

\section{Fioresi's commutation relations}

In \cite{F}, Fioresi has developed useful commutation relations for the
$m\times m$ quantum minors which generate $\gqmn$.  However, Fioresi
works in the following setting.  The field $k$ that she considers is
required to be algebraically closed of characteristic zero, and the
quantum matrix algebra that she considers is generated as an algebra
over the ring $k[q, q^{-1}]$, where $q$ is transcendental over $k$.  The
first thing that we need to do is to observe that these commutation
relations hold over any field $k$ and for any $0\neq q \in k$.  A couple
of warnings about notation for readers comparing \cite{F} with this
paper.  First, because of the choice of relations for $\oqmn$, it is
necessary to replace $q$ by $q^{-1}$ in any relation taken from
\cite{F}.  Secondly, Fioresi works with the quantum grassmannian defined
by the maximal $m\times m$ minors of $\oqnm$; thus, in any maximal
minor, she uses all of the $m$ columns, and a generating quantum minor
of the Grassmannian is specified by choosing $m$ rows.  To deal with
this second difference, we can think of both versions of the quantum
Grassmannian as being subalgebras in the quantum matrix algebra $\oqnn$
and observe that the transpose automorphism, ${\bf \tau}$, see
\cite[3.7.1]{PW}, transforms Fioresi's quantum grassmannian to our
quantum grassmannian.

Recall the following total {\em lexicographic ordering} on   
quantum minors:
$\qs{j_1 j_2 \ldots j_m} \lex \qs{i_1 i_2 \ldots i_m}$ if and only if
there exists an index $\alpha$ such that $j_l=i_l$ for $l<\alpha$, 
but $j_\alpha <i_\alpha$. 

Let $[I] = [i_1 \dots i_m]$ denote an $m\times m$ quantum minor. If
$[I] \neq [1 \dots m]$, consider the least integer $s$ such that $i_s
> s$. Let $\sigma([I])$ be the quantum minor obtained from $[I]$ by
replacing $i_s$ by $i_s -1$ and leaving the other indices unchanged.
Obviously, $\sigma([I]) \lex [I]$. The {\em standard tower} of $[I]$ is
the sequence of quantum minors $[I_N] \blex [I_{N-1}] \blex \dots \blex 
[I_1] 
\blex [I_0]$ where $[I_N] = [I]$, $[I_{l-1}] = \sigma([I_l])$, and $[I_0] =
[1, \dots, r]$. If $[I] = [1 \dots r]$ then the standard tower is
defined to be the single quantum minor $[I]$.

We will denote the version of the $m\times n$ quantum
Grassmannian constructed by Fioresi by $\ghmn$. 
Note also that the relations in \cite{F} use $h$ where we
would use $h^{-1}$; thus we should interchange $h$ and $h^{-1}$.

\begin{proposition}\label{fcomm}
Let $K$ be an algebraically closed field of characteristic zero, and let
$h$ be an indeterminate over $K$. Set $\ghmn$ to be the quantum
Grassmannian subalgebra of $\o_h(M (K[h,h^{-1}])_{mn})$. 
Let $I,J \subseteq \{ 1, \ldots , n \}$ with $|I|=|J|=m$, and
$\qs{I}\lex \qs{J}$. 
Set $s = m- |I\cap J|$. 
Then in $\ghmn$,
\[
\qs{I}\qs{J} = h^s\qs{J}\qs{I} + \sum_{\qs{L}\lex\qs{I}}
\lambda_{\qs{L}}\l h-h^{-1}\r^{i_{\qs{L}}} \l -h\r^{j_{\qs{L}}}     
\qs{L}\qs{L'}, \]
where  $i_{\qs{L}}, j_{\qs{L}} \in \mn$ and $\lambda_{\qs{L}}$ is either
$0$ or $1$,  
while $L'$ is the set $\l I\cap J\r \cup \l ( I\cup J
)\setminus L \r $.
\end{proposition}

\begin{proof} In \cite[Proposition 2.21 and Theorem 3.6]{F}, Fioresi
obtains commutation relations of the above form, but with the products
$\qs{L}\qs{L'}$ on the right hand side of the equation above more
carefully stated.  In Proposition 2.21 she first obtains the result for
the case that $I\cap J = \emptyset$.  In this case, the quantum minors
$[L]$ involved are members of the standard tower of $[I]$, and so
$\qs{L}\lex\qs{I}$, as we require.  The general case where $I\cap J \neq
\emptyset$ is dealt with in Theorem 3.6.  Set $[\widetilde{I}]$ to be the
quantum minor obtained from columns $I\backslash 
(I\cap J)$, and similarly, define
$[\widetilde{J}]$.  Proposition 2.21 provides a commutation rule for
$[\widetilde{I}][ \widetilde{J}]$ with terms on the right hand side
$[\widetilde{L}][\widetilde{L'}]$ where $[\widetilde{L}] \lex
[\widetilde{I}]$. In Theorem 3.6, a commutation rule with the same
coefficients is then obtained
for $\qs{I}\qs{J}$ by replacing each $[\widetilde{L}][\widetilde{L'}]$
by $[\widetilde{L}\cup (I\cap J)][\widetilde{L'}\cup (I\cap J)]$. Thus,
all that needs to be done is to make the easy observation that if
$[\widetilde{L}] \lex
[\widetilde{I}]$ then $[\widetilde{L}\cup (I\cap J)] \lex
[\widetilde{I}\cup (I\cap J)] = [I]$. 
\hspace{2ex}\smallbox\end{proof} 

\begin{corollary}\label{commr}
Let $k$ be any field and $q$ any nonzero element of $k$. 
Set $\gqmn$ to be the quantum
Grassmannian subalgebra of $\oqmn$. 
Let $I,J \subseteq \{ 1, \ldots , n \}$ with $|I|=|J|=m$, and  
$\qs{I}\lex \qs{J}$. Set $s= m - |I\cap J|$. Then in $\gqmn$,
\[
\qs{I}\qs{J} = q^s\qs{J}\qs{I} + \sum_{\qs{L}\lex\qs{I}}
\lambda_{\qs{L}} \l q-q^{-1}\r^{i_{\qs{L}}} \l -q\r^{j_{\qs{L}}}
\qs{L}\qs{L'}, \] 
where $\lambda_{\qs{L}} \in k$, $i_{\qs{L}}, j_{\qs{L}} \in \mn$ 
and $\lambda_{\qs{L}}$ is either
$0$ or $1$,
while $L'$ is the set $\l I\cap J\r \cup \l ( I\cup J
)\setminus L \r $.
\end{corollary}

\begin{proof} Proposition~\ref{fcomm} applies in the case that $K =
\mc$. In this case, 
observe that the coefficients of the monomials in the maximal minors 
are all in $\mz[h,h^{-1}]$; so that these
relations hold in the quantum Grassmannian over $\mz[h,h^{-1}]$.  There
is then a natural homo\-morphism from this quantum Grassmannian to
$\gqmn$, such that $z\mapsto z1_k$ for $z \in \mz$ 
and $h\mapsto q$, which produces the
required relations. 
\hspace{2ex}\smallbox\end{proof} \\

Recall that an element $a$ of an algebra $A$ is a {\em normal} element
if $aA = Aa$. 
The next result follows immediately from the previous Corollary.

\begin{corollary}\label{daft}
An $m\times m$ 
quantum minor $\qs{I} \in \gqmn$ is normal modulo the ideal generated by
the set $\{ \qs{J} \mid \qs{J}\lex\qs{I}\}$.
\end{corollary} 

The algebra $\oqmn$ is a connected $\mn$-graded
algebra, graded by the
total degree in the canonical generators. Since $\gqmn$ is a subalgebra
generated by homogeneous elements of degree $m$ with respect to this
grading, $\gqmn$ inherits a connected $\mn$-graded structure in which
its canonical generators have degree one. 

\begin{theorem} The quantum grassmannian $\gqmn$ is a noetherian domain.
\end{theorem}

\begin{proof} The quantum Grassmannian $\gqmn$ is generated by the
${{n}\choose{m}}$ quantum minors of size $m$ in $\oqmn$.  Denote these
quantum minors
by
$u_1 \lex u_2 \lex \ldots \lex
u_{n\choose m}$. Then by Corollary~\ref{daft}, $\{ u_1,
\ldots ,u_{{n}\choose{m}}\}$ is a normalising sequence of $\gqmn$;
that
is, $u_1$ is normal and $u_l$ is normal modulo 
the ideal generated by $\{ u_1, \ldots ,
u_{l-1}\}$, for $l>1$. 
The factor by the ideal generated by this normalising sequence is the
base field; so the fact that $\gqmn$ is noetherian follows by repeated
use of \cite[Lemma 8.2]{ATV}. 

Finally, $\gqmn$ is a domain since it is a  subalgebra of $\oqmn$ which is a
domain. 
\hspace{2ex}\smallbox\end{proof} \\

\noindent{\bf Remark}~~
If $A$ is a noetherian, connected
$\mathbb{N}$-graded $k$-algebra such that every non-simple graded prime
factor ring $A/P$ contains a nonzero homogeneous normal element in
$\oplus_{i\geq 1} \l A/P\r_i$ then we say that $A$ has
{\em enough normal elements} (\cite{Z}).
Thus, the two previous results show that the
quantum grassmannian has enough normal elements.\\

There is a useful isomorphism between $\gqmn$ and $\gqminusmn$ which we
now describe.  
Notice that, if $1 \le i_1 < \dots < i_m \le n$, ${\cal G}_q(m,n)$ is
isomorphic to the subalgebra of ${\cal O}_q(M_n)$ generated by the
$m \times m$ minors that use rows $i_1 , \dots , i_m$, that is, the
minors $[I|J]$ with $I=\{i_1 , \dots , i_m\}$ and $J \subseteq
\{1,\dots,n\}$, $|J|=m$.
Let $A:=
\oqnn$ with generators $X_{ij}$ and $A':= \oqminusnn$ with generators
$X_{ij}'$.  Take a copy $R$ of $\gqmn$ inside $A$ generated by the
$m\times m$ quantum minors that use the first $m$ rows of $A$, and take
a copy $R'$ of $\gqminusmn$ that uses the last $m$ rows of $A'$. 
Following the proof of \cite[Corollary 5.9]{GL2}, we see that there is
an isomorphism $\bfdelta: A\goesto A'$ which takes $[I|J]$ to
$[\omega_0I|\omega_0J]'$, where $[-|-]'$ denotes a quantum minor in
$A':= \oqminusnn$ and $\omega_0$ is the longest element of the symmetric
group $S_n$; that is, $\omega_0(i) = n-i+1$.  Note that the isomorphism
$\bfdelta$ restricted to $R$ produces an isomorphism from $R$ to $R'$
that takes a generating minor $[I]$ to the minor $[\omega_0I]'$.  In
particular, note that under this isomorphism, $[1 2  \dots m]$, the
leftmost minor of $R = \gqmn$, is translated into the rightmost minor
$[n-m+1 \dots n]'$ of the quantum grassmannian $R'= \gqminusmn$.  We
denote this induced isomorphism from $\gqmn$ to $\gqminusmn$ by
$\bfdelta$ also. 

As an example of the use of the isomorphism $\bfdelta$, we record the
following lemma which we need later.

\begin{lemma}\label{mincomm}
Let $I\subseteq \{ 1, \ldots ,n \}$ with $|I|=m$.  Then
\[
\qs{I}[n-m+1 \dots n] = q^{s}[n-m+1 \dots n]\qs{I}
\]
where $s= m -| I \cap \{ n-m+1, \ldots ,n\}|$, and thus $[n-m+1 \dots n]$ 
is normal in
$\gqmn$.
\end{lemma}

\begin{proof} 
Note that $\omega_0\{n-m+1, \dots, n\} = \{1, \dots, m\}$. Note also
that \newline 
$| I \cap \{ n-m+1, \ldots ,n\}| = |\omega_0I \cap \omega_0\{n-m+1,
\dots, n\}| = |\omega_0I \cap \{1, \dots, m\}|$. 

By Corollary~\ref{commr}, 
$[1 \dots m][\omega_0I] = q^s[\omega_0I][1 \dots m]$. 
Applying $\delta$ to this equation gives
$[n-m+1 \dots n]'[I]' = q^s[I]'[n-m+1 \dots n]' $ 
in $\gqminusmn$. This can be rewritten as 
$[I]'[n-m+1 \dots n]' = q^{-s}[n-m+1 \dots n]'[I]'$ 
in $\gqminusmn$. 
Finally, replacing $q^{-1}$ by $q$, we obtain 
\[
[I][n-m+1 \dots n] = q^{s}[n-m+1 \dots n][I]
\]
in $\gqmn$.
\hspace{2ex}\smallbox\end{proof}

\section{A basis for $\gqmn$}

In this section, we obtain a basis for $\gqmn$.  This basis is a subset
of the basis of  preferred products of $\oqmn$ obtained in
\cite[Section 1]{GL}.  First, we adapt the language used in that paper to the
grassmannian subalgebra $\gqmn$.  Recall from Section 1 that if $J$ is
an $m$-element subset of $\{1, \dots, n\}$ then $[J]$ denotes the
quantum minor $[1, \dots, m \mid J]$ of $\oqmn$.  Thus, 
let $m,n \in \mn^{\ast}$ with $n \geq m$.  We define a partial ordering
on $m$-element subsets of $\{ 1, \ldots ,n\}$.

\begin{definition} 
 Let
$A,B\subseteq \{1,\dots,n\}$ with $|A| = m = |B| $.
We define a partial ordering, denoted by $\le_*$. Write $A$
and $B$ in ascending order:
$$A= \{a_1< a_2< \dots <a_m\} \qquad \text{and} \qquad B= \{b_1< 
b_2< \dots <b_m\}.$$
Define $A \le_* B$ to mean that $a_i\le b_i$ for 
$i=1,\dots,m$.
\label{po}
\end{definition} 

This naturally defines a partial ordering on the generators of $\gqmn$.
\begin{definition}
Let $\left[ I\right]$ and $\left[ J\right]$ belong to the generating
set of $\gqmn$.  Then we write that $\qs{I} \leq_c \qs{J}$ if and only
if $I\le_* J$.
\end{definition} 

\begin{figure}[ht]
\ignore{
$$\xymatrixrowsep{2.4pc}\xymatrixcolsep{3.2pc}\def\objectstyle{\scriptstyle}
\xymatrix@!0{
 && [ 456 ] \edge[d]\\
 && [ 356] \edge[dl] \edge[dr]\\
 & [ 346] \edge[dl] \edge[dr] &&[ 256] \edge[dl] \edge[dr]\\
[ 345] \edge[dr] && [ 246] \edge[dl] \edge[dr] \edge[d] && [ 156] 
  \edge[dl]\\
 & [ 245] \edge[d] \edge[dr] & [ 236] \edge[dl]|\hole \edge[dr]|\hole &[146]
\edge[dl] \edge[d]\\
 &[235] \edge[dl] \edge[dr] &[145] \edge[d] &[136] \edge[dl] \edge[dr]\\
[234] \edge[dr] &&[135] \edge[dl] \edge[dr] &&[126]\edge[dl]\\
 &[134] \edge[dr] &&[125] \edge[dl]\\
 &&[124] \edge[d]\\
 &&[123]
}$$
}
\caption{The partial ordering $\leq_c$ on $\gq(3,6)$
}\label{gra}
\end{figure}
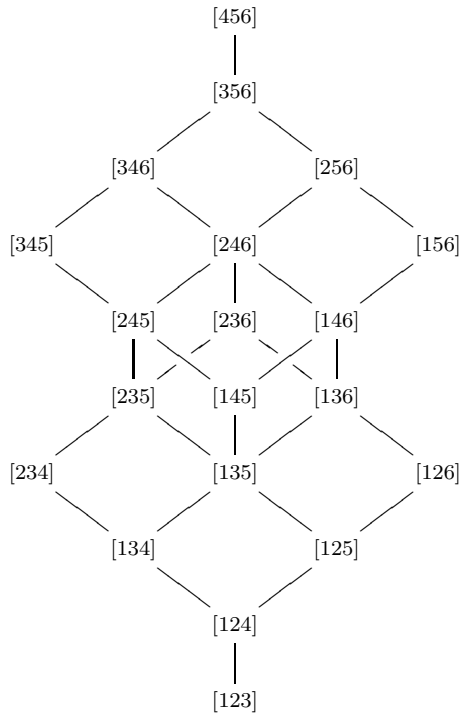

For example, Figure~\ref{gra} shows the ordering on generators of
$\gq(3,6)$.\\ 

Recall that a {\em tableau} is a Young diagram with entries in each box. 
If each row of a tableau $T$ has length $m$ then we will say that $T$ is
an {\em $m$-tableau}.  Here, we consider tableaux with entries from $\{1,
\dots, n\}$ and no repetitions in each row.  An {\em allowable
$m$-tableau} $T$ is an $m$-tableau with strictly increasing rows.  If an
allowable $m$-tableau $T$ has rows $J_1, \dots, J_s$, then $T$ is {\em
preferred} if and only if $J_1 \leq_* J_2 \leq_* \ldots \leq_* J_s$.

Let $I = \{m, m-1, \dots, 1\}$ and let $S$ be an $m$-tableau which has
the same number of rows as $T$ and such that each row of $S$ is $I$. Then
$T$ is  an allowable (preferred) $m$-tableau if and only if the
bitableau $(S\mid T)$ is allowable (preferred) in the sense of
\cite{GL}. With this in mind, we define the following ordering on
allowable $m$-tableau. Let 
\[
T = 
\begin{pmatrix}
J_1\\
J_2\\
\vdots\\
J_t
\end{pmatrix}, \qquad 
S =
\begin{pmatrix}
L_1\\
L_2\\
\vdots\\
L_s
\end{pmatrix}.
\]
Then $T \prec S$ if $t>s$, or if $s=t$ and
\[
\{ J_1, \ldots ,J_t \} \lex \{ L_1, \ldots ,L_s \};
\]
that is, there exists an index $i$ such that $J_\alpha =
L_\alpha$ for $\alpha < i$, but $J_i <_* L_i$.

Any allowable $m$-tableau determines a product of quantum minors in
the quantum grassmannian as follows.
\begin{definition} 
For any (allowable) $m$-tableau
\[ 
T =
\begin{pmatrix}
J_1\\
J_2\\
\vdots\\
J_s
\end{pmatrix},
\]
define  $[T] = [J_1][J_2]\ldots [J_s]$.
\end{definition} 
\begin{definition} 
The {\em content} of an $m$-tableau $T$ is the multiset 
$\{1^{t_1}, 2^{t_2}, \cdots ,n^{t_n}\}$, where $t_i$ is the
number of times $i$ appears in $T$. 
\end{definition} 

We will use the content of a tableau to define a natural
$\mathbb{Z}^n$-grading on the $m\times n$ quantum Grassmannian. 
There is a $\mathbb{Z}^n$-grading on $\oqmn$ defined by assigning degree
$\varepsilon_j$ to $X_{ij}$, where $\varepsilon_j$ for $j=1, \dots, n$
form the natural basis of $\mathbb{Z}^n$. Since the maximal minors of
$\oqmn$ are homogeneous with respect to this basis, there is an induced
$\mathbb{Z}^n$-grading on $\gqmn$: consider a product of minors $\qs{T}$
in $\gqmn$, if the tableau
$T$ has content 
$\{ 1^{t_1}, 2^{t_2}, \ldots ,n^{t_n}\}$, then 
$\qs{T}$
is homogeneous of degree  $\l t_1, t_2, \ldots ,t_n\r$.  
Thus, the degree of a product is dependent on the number of times each 
column of the $m\times n$ quantum matrix appears in it.

\begin{theorem}\label{gqpr}
Generalised Quantum Pl\"{u}cker Relations for Quantum Grassmannians\\ 
Let  $J_1,J_2,K \subseteq \{ 1, 2, \ldots
,  n \} $ be such that $|J_1|, |J_2| \leq m$ and 
$|K| = 2m-|J_1| - |J_2| > m$.  Then
\begin{equation*}
\sum_{K' \sqcup K'' = K }^{} \l -q \r ^{ \ell \l J_1;K'\r  +\ell\l K';K''\r 
+\ell\l K'';J_2\r } [J_1 \sqcup K' ] [K'' \sqcup J_2 ] = 0,
\end{equation*}
where
$ \ell \l I;J\r = | \{ \l i,j\r \in I \times J : i >j \}|$.
\end{theorem}

\begin{proof}
We work in the algebra 
$\oqnn$ and apply \cite[Proposition B2(a)]{GL} with $I_1 =
I_2 = \{1, \dots, m\} =:I$. Thus, 
\begin{equation*}
\sum_{K' \sqcup K'' = K }^{} \l -q \r ^{ \ell \l J_1;K'\r  +\ell\l
K';K''\r
+\ell\l K'';J_2\r } [I|J_1 \sqcup K' ] [I|K'' \sqcup J_2 ] = 0,
\end{equation*}
since $|K| > m = |I_1\cup I_2|$, see \cite[B3]{GL}. This is the desired
relation. 
\hspace{2ex}\smallbox\end{proof}

\begin{lemma}\label{span1}
 Let $T $ be an $m$-tableau with content
$\gamma$ and suppose that $T$ is not preferred.  Then

{\rm (a)} $T$ is not minimal with respect to $\prec$ among
$m$-tableaux with content $\gamma$;

{\rm (b)} $[T]$ can be expressed as a linear combination of
products
$[S]$, where each $S $ is an $m$-tableau with content
$\gamma$ such that
$S \prec T$.
\end{lemma}

\begin{proof}
Follow the proof of \cite[Lemma 1.7]{GL}. Note that in the proof the
only place where the shape of a bitableau might change is near the end
of the proof where the right-hand side of the Exchange Formula is
considered. In our situation, the right-hand side is zero, as noted in
Theorem~\ref{gqpr}. 
\hspace{2ex}\smallbox\end{proof}\\

Note that fixing the content of an $m$-tableau fixes its shape and thus
fixes the number of rows in the $m$-tableau. 

Let $\delta = (c_1, \dots, c_n) \in \mn^n$. Let $V$ be the homogeneous
component of degree $\delta$ in $\gqmn$. Note that $V$ might be zero,
and that this is the case if and only if there is no product $[T]$ where
$T$ is an $m$-tableau of content $(1^{c_1}\dots n^{c_n})$. Also, an
element of $\gqmn$ belongs to $V$ if and only if it is a linear
combination of products $[T]$, where $T$ runs over all $m$-tableau
with content $(1^{c_1}\dots n^{c_n})$; that is, the  products $[T]$,
where $T$ runs over all $m$-tableau
with content $(1^{c_1}\dots n^{c_n})$ span $V$. 

\begin{theorem}\label{grbasis}
Let $\delta = \l  c_1,\ldots,c_n \r \in \mn^n$, 
let $V$ be the homogeneous
component of $\gqmn$ with degree $\delta$, and set 
 $\gamma= \l 1^{c_1} 2^{c_2} \cdots n^{c_n}\r$.  The products $[T]$, as
$T$ runs over all preferred $m$-tableau with content $\gamma$, form a
basis for $V$.
\end{theorem}

\begin{proof}
It is enough to prove that for any $m$-tableau $T$ 
with content $\gamma$ the
product $[T]$ is a linear combination of products $[S]$ where $S$ is a
preferred $m$-tableau with content $\gamma$. Let $\ce$ be the set of
$m$-tableau with content $\gamma$; clearly, $\ce$ is a finite set and we
order it by $\prec$. We use induction on $\prec$ to show the result. Let
$U\in \ce$. If $U$ is minimal, then it is preferred, by part (a) of 
the previous result. Otherwise, by part (b) of the previous result,
$[U]$ is a linear combination of products $[S]$, where $S\in\ce$ and
$S \prec U$. Thus, by an induction argument applied to $S$, we may
conclude that $[U]$ is a linear combination of products $[S]$ where
$S$ is a preferred $m$-tableau with content $\gamma$.
 
Recall that $\gqmn$ is a subalgebra of $\oqmn$ and notice that 
the products $[T]$, as $T$ runs over all preferred $m$-tableaux of
content $\gamma$, form a
subset of the basis of $\oqmn$ constructed in \cite{GL}.  Therefore,
they are linearly independent and we have the result. 
\hspace{2ex}\smallbox\end{proof}

\begin{corollary}\label{gbasis}
 The products $[T]$, as
$T$ runs over all preferred $m$-tableaux, form a basis for
$\gqmn$. \end{corollary}

This basis can be used to calculate the Gelfand-Kirillov dimension of
the $m\times n$ quantum Grassmannian.

Consider the partial ordering $\leq_c$ on the generating 
minors of
$\gqmn$.  A {\em saturated path} between two
minors $a<_c b$ will be an `upwards path' $a=a_1<_c a_2 <_c \ldots <_c a_l =b$
of minors such that no additional terms can be added; that is, for any index
 $i$ there is no minor $d$ such that $a_i<_c d<_c a_{i+1}$.  The
{\em length} of such a saturated path is defined to be $l$.
For example, a saturated path between the minors $\qs{134}$ and
$\qs{256}$ in $\gq(3,6)$ is
\[
\qs{134} <_c \qs{234} <_c \qs{235} <_c \qs{236} <_c
\qs{246} <_c \qs{256}.
\]
The length of this saturated path is $6$.

A {\em maximal path} is a saturated path between the two minors
$\qs{1\ldots m}$ and $\qs{n-m+1 \ldots n}$.  It is easy to check that
any maximal path has length  $m(n-m)+1$.

\begin{proposition}\label{gragk} 

Let $G= \gqmn$ and let $\alpha$ be the length of a maximal path in $G$. 
Then 
\[ \gkdim(\gqmn) =\alpha = m(n-m)+1 .  
\]
\end{proposition} 

\begin{proof}  
Let $V$ be the $k$-subspace of $G$
spanned by the $m\times m$ minors which generate $G$. Then 
$\gkdim(G) =
\overline{\lim}\log_n d_V\l n\r$ where 
$d_V\l n\r =
\dim_k\l \sum_{i=0}^{n} V^i\r$. 
Let $a_1, a_2, \ldots , a_\alpha$ be a
maximal path in $G$.  Then $a_1^{s_1}a_2^{s_2}\ldots
a_\alpha^{s_\alpha} \in V^{n+1}$ whenever $\sum_{i=1}^{\alpha} s_i =
n+1$. The set $\{ a_1^{s_1}a_2^{s_2}\ldots a_\alpha^{s_\alpha} 
\mid \Sigma
s_i = n+1 \}$ is linearly independent.  Therefore
\begin{align*} \dim_k\l
V^{n+1}\r \geq& \;|\{ a_1^{s_1}a_2^{s_2}\ldots a_\alpha^{s_\alpha} 
\mid \Sigma
s_i = n+1 \}|= \;{{n+\alpha} \choose {\alpha -1}} 
\end{align*}
which is a polynomial in $n$ of degree $\alpha -1$.  It follows that
$\gkdim(G) \geq \alpha$. 

Let $a_{i_1} \dots a_{i_n} \in V^n$.  
By Theorem~\ref{grbasis}, 
$a_{i_1} \dots a_{i_n}$ may be rewritten as a linear
combination of preferred products from $V^n$.

There are finitely many maximal paths in $\gqmn$.  Suppose there are $c$
such paths and index them  $1,\ldots ,c$.  Let $a_1<_c a_2 <_c\ldots
<_c a_\alpha$ be the $i$th maximal path and let $W_i^{(n)}$ denote the subspace
generated by monomials
$a_1^{s_1}\ldots a_\alpha^{s_\alpha}$ such that $ \Sigma
s_j =n$. The above observation shows that 
$V^n \subseteq \sum_{i=1}^{c} W_i^{(n)}$. Consider ${\rm{dim}}(W_i^{(n)})$.
The products $a_1^{s_1}\ldots a_\alpha^{s_\alpha}$ such that
$\Sigma s_j =n$ are linearly independent.  Therefore
\begin{eqnarray*}
{\rm{dim}}(W_i^{(n)}) &=& |\{  a_1^{s_1}\ldots a_\alpha^{s_\alpha}\mid \sum s_i
=n\}|
= |\{ (s_1, \ldots ,s_\alpha )\in \mn^\alpha \mid \sum s_i=n\}|.
\end{eqnarray*}
Therefore   ${\rm{dim}}(W_i^{(n)})=  {\rm{dim}}(W_j^{(n)})$ 
for all $i,j\in \{1,
\ldots ,c\}$.   Thus
\[
{\rm{dim}}\l V^n\r \leq {\rm{dim}}\l \sum_{i=1}^{c}W_i^{(n)} \r
\leq c\, {\rm{dim}} \l W_1^{(n)}\r = c\,{{n+\alpha-1} \choose {\alpha -1}}
\]
and $d_V\l n\r \leq c\sum_{i=0}^n {{i+\alpha -1} \choose {\alpha -1}}$, 
a polynomial of degree $\alpha$. 
It follows that $\gkdim(G) \leq \alpha$.  Hence, $\gkdim(G) =
\alpha = m(n-m) +1$. 
\hspace{2ex}\smallbox\end{proof}\\ 

For example, $\gkdim(\gqtwofour) = 2(4-2) +1 = 5$.

\section{Noncommutative Dehomogenisation} 

If $R$ is a commutative $\mn$-graded algebra, and $x$ is a homogeneous
nonzerodivisor in degree one, then the {\em dehomogenisation} of $R$ at
$x$ is usually defined to be the factor algebra $R/(x-1)R$,
\cite[Appendix 16.D]{BV}.  This
definition is unsuitable in a noncommutative algebra if the element $x$
is merely normal rather than central: in this case, the factor algebra
is often too small to be useful.  For example, let $R$ be the quantum
plane $k_q[x,y]$ with $xy = qyx$ and $q\neq 1$.  Setting $x=1$ forces
$y=0$; so that the factor algebra $R/\langle x-1\rangle$ 
is isomorphic to $k$ rather
than being a one-dimensional algebra, as one might hope.  However, in
the commutative case, an alternative approach is to observe that the
localised algebra $S:=R[x^{-1}]$ is $\mz$-graded, $S = \oplus_{i\in \mz}\,
S_i$, and that $S_0 \cong R/(x-1)R$.  If $x$ is a normal nonzerodivisor
of degree one 
in a noncommutative $\mn$-graded algebra $R = \oplus_{i\in \mn}\,
R_i$, then one can form the Ore localisation
$R[x^{-1}] =:S$, and then this second approach does yield a
useful algebra in the noncommutative case. Indeed, for $i,j\in\mn$
denote by $R_ix^{-j}$ the $k$-subspace of elements of $S$ that can be
written as $rx^{-j}$ with $r\in R_i$; clearly, $R_ix^{-j} \sse
R_{i+1}x^{-(j+1)}$. For $l\in \mz$, set $S_l = \sum_{t\geq 0}\,
R_{l+t}x^{-t} = \cup_{t\geq 0}\, R_{l+t}x^{-t}$. Then $S$ is a
$\mz$-graded algebra with $S=\oplus_{l\in\mz}\, S_l$. 

\begin{definition} Let $R = \oplus R_i$ be an $\mn$-graded $k$-algebra 
and let $x$ be
a regular homogeneous normal element of $R$ of degree one.  Then the {\em
dehomogenisation} of $R$ at $x$, written $\dhom(R,x)$, is defined to be
the zero degree subalgebra $S_0$ of the $\mz$-graded algebra $S:=
R[x^{-1}]$. 
\end{definition}

It is easy to check that $\dhom(R,x) = \sum_{i=0}^{\infty}\, R_ix^{-i} 
= \cup_{i=0}^{\infty}\, R_ix^{-i}$.
In particular, if $R = k[R_1]$ then $\dhom(R,x) = \sum_{i=0}^{\infty}\,
(R_1x^{-1})^i = \cup_{i=0}^{\infty}\,
(R_1x^{-1})^i$, and further, if $R_1 = ka_1 + \dots + ka_s$ then
$\dhom(R,x) = k[a_1x^{-1}, \dots, a_sx^{-1}]$. 

Denote by $\sigma$ the automorphism of $S$ given by $\sigma(s) =
xsx^{-1}$ for $s\in S$. Note that $\sigma$ induces an automorphism of
$S_0$, also denoted by $\sigma$. 

\begin{lemma} \label{sle}
Let $R$ be an $\mn$-graded algebra  and let $x$ be a regular normal
homogeneous
element of degree $1$.  Then there is an isomorphism 
$$\theta:\dhom(R,x)[y,y^{-1}; \sigma] \goesto R[x^{-1}]$$ 
which is the identity on
$\dhom(R,x)$ and sends $y$ to $x$. 
\end{lemma} 

\begin{proof} The existence of $\theta$ is clear from the universal
property of skew-Laurent extensions. It is easy to check that $\theta$
is an isomorphism. 
\hspace{2ex}\smallbox\end{proof} \\

Some properties of dehomogenisation follow in an elementary way from
this result.

\begin{corollary} \label{dhomprops}
Let $R = \oplus_{i\geq 0}\, R_i$ 
be an $\mn$-graded algebra and let $x$ be a
regular homogeneous normal element of degree one.  
\newline (i) $R$ is a domain if and only if $\dhomrx$ is a domain. 
\newline (ii) If $R$ is noetherian then $\dhomrx$ is noetherian. 
\newline (iii) If $R$ is locally finite (that is, $\dim(R_i) < \infty$
for all $i\in \mn$) then 
$\gkdim(R) = \gkdim(\dhomrx) +1$. 
\end{corollary} 

\begin{proof}
Point (i) follows at once from the isomorphism in Lemma~\ref{sle}. 

\noindent(ii)~~If $R$ is noetherian then so is $R[x^{-1}]$ and thus
$\dhom(R,x)[y,y^{-1}; \sigma]$ is noetherian, by Lemma~\ref{sle}. As is
well-known, since $\sigma$ is an automorphism of $\dhom(R,x)$, this
implies that $\dhom(R,x)$ is noetherian. 

\noindent(iii)~~Let $\sigma $ be the automorphism of $R$ induced by
conjugation by $x$. It is clear that $\sigma$ is a graded automorphism; 
and so from the local finiteness of $R$, we see that the elements 
$x^i$, for $i\geq
1$, are  local normal elements in the sense of \cite[p168]{KL}. 
By using \cite[12.4.4]{KL}, it follows
that $\gkdim(R[x^{-1}]) = \gkdim (R)$. On the other hand, the automorphism
$\sigma$ induced on $S_0$ by conjugation by $x$ in $S$ is locally
algebraic in the sense of \cite[p164]{KL}. 
Indeed, $S_0 = \cup_{t\geq 0}\, R_tx^{-t}$ and for all
$t\in\mn$ the $k$-subspace $R_tx^{-t}$ is a finite dimensional
$\sigma$-stable subspace of $S_0$. It follows from \cite[p164]{KL} that
$\gkdim (S_0[y,y^{-1}; \sigma]) = \gkdim (S_0) + 1$. The conclusion
follows from Lemma~\ref{sle}. 
\hspace{2ex}\smallbox\end{proof}

\section{Dehomogenisation of $\gqmn$}\label{dehg}

In the classical commutative theory it is a well-known and basic result
that the dehomogenisation of the homogeneous 
coordinate ring of the $m\times n$
Grassmannian at the minor $[n-m+1, \dots, n]$ is isomorphic to the
coordinate ring of $m\times (n-m)$ matrices; that is,
\[
\frac{\o(\cg(m,n))}{\left< [n-m+1, \dots, n] -1\right>} \cong
\o(M_{m,n-m}(k)).
\]

In this section, we show that the corresponding result holds for
$\gqmn$ when we use the noncommutative dehomogenisation defined in
the previous section. Recall from Lemma~\ref{mincomm} 
that $[n-m+1, \dots, n]$ is
a normal element of $\gqmn$: in fact, it $q$-commutes with the other
maximal minors, and this will be important in calculations.  

Recall that we may 
consider $\gqmn$ to be a $\mn$-graded algebra with each $m\times
m$ quantum minor given degree $1$. Set $x = [n-m+1, \dots, n]$ and  
$S:= \gqmn[x^{-1}]$, and note
that $\dhom(\gqmn, [n-m+1, \dots, n]) = S_0$ is generated by elements
of the form $\{ I\} := \qs{I}[n-m+1, \dots, n]^{-1}$ with $I \sse \{1,
\dots, n\}$ and $|I|=m$, see Section 3. \\

Now let $u$ be a positive integer and consider $\oquu$.
If $I\subseteq \{1, \dots, u\}$ then $\Itil:= \{1, \dots, u\}\backslash
I$. In an exponent $I$ denotes the sum of the indices occuring in the
index set $I$. 

Let $D_q$ be the quantum determinant of
$\oquu$. Since $D_q$ is a central element, we can invert it to form the
{\em $u\times u$ quantum general linear group} $\oqglu :=
\oquu[D_q^{-1}]$.  The algebra $\oqglu$ is a Hopf algebra, with antipode
$S$, and counit $\varepsilon$. 

There is a useful antiendomorphism $\Gamma: \oquu \goesto \oquu$ defined
on generators by $\Gamma(X_{ij}) = (-q)^{i-j}[\widetilde{\{j\}} |
\widetilde{\{i\}}]$, see \cite[Corollary 5.2.2]{PW}. We need to know the
effect of $\Gamma$ on quantum minors. This is given in the following
lemma, which is presumably well-known but we give a proof since we have
been unable to find a clear exposition. Recall that $\Delta([I | J]) =
\sum_{|K| = |I|}\, [I|K]\ot[K|J]$, where $\Delta$ is the
comultiplication map on $\oquu$, by \cite[(1.9)]{NYM}. Recall also that
$\varepsilon([I|J])$ equals $1$ if $I=J$ and
$0$ otherwise.

\begin{lemma} \label{Gamma}
 Let $[I|J]$ be an $r\times r$ quantum minor in $\oquu$. Then,

\noindent (i)~~$S([I|J]) = (-q)^{I-J}[\Jtil|\Itil]D_q^{-1}$

\noindent (ii)~~$\Gamma([I|J]) = (-q)^{I-J}[\Jtil|\Itil]D_q^{r-1}$

\end{lemma} 

\begin{proof} 
We establish the  first claim by
calculating the expression
$$
\sum_{K,L}\, (-q)^{L - J}S([I|K])[K|L][\Jtil|\Ltil]D_q^{-1}
$$ 
in two different ways.

First,
\begin{eqnarray*}
\sum_{K,L}\, (-q)^{L - J}S([I|K])[K|L][\Jtil|\Ltil]D_q^{-1}
 &=&
\sum_{K}\, S([I|K])\left\{\sum_{L}\,(-q)^{L -
J}[K|L][\Jtil|\Ltil]D_q^{-1}\right\} \\
 &=&
\sum_{K}\, S([I|K])\varepsilon([K|J])1 = S([I|J]),
 \end{eqnarray*}
by using the first equality of \cite[4.4.3]{PW}.

Secondly,
\begin{eqnarray*}
\sum_{K,L}\, (-q)^{L - J}S([I|K])[K|L][\Jtil|\Ltil]D_q^{-1}
        &=&
\sum_{L}\, \left\{ \sum_{K}\,S([I|K])[K|L] \right\} (-q)^{L -
J}[\Jtil|\Ltil]D_q^{-1} \\
        &=&
\sum_{L}\, \varepsilon([I|L]) (-q)^{L -
J}[\Jtil|\Ltil]D_q^{-1} \\
        &=&  (-q)^{I -
J}[\Jtil|\Itil]D_q^{-1},
 \end{eqnarray*}
by using the defining property of the antipode.

The second claim follows easily from the first, since $S([I|J]) =
\Gamma([I|J])D_q^{-r}$ for $r\times r$ quantum minors $[I|J]$.  This is
easily established from  the fact that it holds on the generators $X_{ij}$
and that $S$ and $\Gamma$ are anti-endomorphisms. 
\hspace{2ex}\smallbox\end{proof} \\

 We will  need
the anti-endomorphism $\Gamma\circ\tau: \oquu \goesto \oquu$ defined by
$\Gamma\circ\tau(X_{ij}) = (-q)^{j-i}[\widetilde{\{i\}} |
\widetilde{\{j\}}]$ for $1\leq i,j\leq u$. Here, $\tau$ is the
transposition automorphism given in \cite[Proposition 3.7.1(1)]{PW}.
Note that, by Lemma~\ref{Gamma},  
the effect of $\Gamma\circ\tau$ on the $r\times r$  
quantum minor $[I\mid J]$ is given by 
$\Gamma\circ\tau([I\mid J]) =
(-q)^{J-I}[\widetilde{I}\mid\widetilde{J}]D_q^{r-1}$. 

Given $I=\{ i_1, \ldots , i_m\}\subseteq \{ 1, \ldots , n\}$  
the set $I\backslash \{i_k\}$ 
is denoted by $\{ i_1, \ldots , \widehat{i_k}, \ldots , i_m\}$. 
Given two sets $I, J \subseteq \{ 1, \ldots , n\}$ recall that 
\[
\ell \l
I;J\r := | \{ \l i,j\r \in I \times J : i >j \}|.
\]
In the next proof, and throughout the paper, $\qdot$ denotes a power of
$-q$ that is not necessary to keep track of explicitly. 

\begin{lemma} \label{gens}
The $k$-algebra $\dhom(\gqmn, [n-m+1, \dots, n]) = S_0$ is
generated as an algebra by the elements $\{j \quad n-m+1 \ldots
\widehat{i} \ldots n\}$ for $1 \leq j \leq n-m <  i \leq n$. 
\end{lemma} 

\begin{proof} 
We know that $S_0$ is generated by the elements 
$\{ I\} := \qs{I} [n-m+1, \dots, n]^{-1}$, 
where $I\sse \{1, \dots, n\}$ and $ |I| =m$.  
We show that each such element can be
expressed as a $k$-linear combination of products of elements of the
form $ \{j \quad n-m+1 \ldots \widehat{i} \ldots n\}$, 
where $ 1 \leq j \leq n-m <  i \leq n$. Denote by $A$ the subalgebra 
of $S_0$ generated  by the elements 
$ \{j \quad n-m+1 \ldots \widehat{i} \ldots
n\}$. 

Let $I = \{ i_1\leq \ldots \leq i_m\} \neq \{n-m+1, \dots, n\}$ 
be an ordered subset of $\{ 1,
\ldots , n\}$ and let $2\leq t\leq m+1$ be such that $i_t\geq n-m+1$ but
$i_{t-1}<n-m+1$; that is, $I\cap \{1, \ldots , n-m\} = \{ i_1, \ldots
i_{t-1}\}$.  We will use induction on $t$ to show that $\{I\} \in A$. 

If $t=2$, then $I$ is of the form $ \{j \quad n-m+1 \ldots \widehat{i}
\ldots n\}$ and so $\{ I\}\in A$.  Consider a fixed $t \in \{3, \dots,
m+1\}$ and suppose that the result is true for $t-1$. 
Now consider $\qs{I}= \qs{i_1\ldots i_m}$ with $I\cap  \{1,
\ldots , n-m\}= \{ i_1, \ldots i_{t-1}\}$.  We use the generalised Quantum
Pl\"{u}cker relations (Theorem \ref{gqpr}) to rewrite the product
$ [n-m+1, \dots, n]\qs{i_1\ldots i_m}$.

Let $K= \{ i_1, n-m+1, \ldots, n\}$, $J_1=\emptyset$ and $J_2=\{ i_2,
\ldots ,i_m\}$.  Then
\begin{align*}
\\[-3ex]
\sum_{K'\sqcup K'' =K} \l -q\r^{\ell\l K';K''\r +\ell\l K'';J_2\r}
\qs{K'}\qs{K''\sqcup J_2} =0
\end{align*}
where either
\[
K'=\{ n-m+1,\ldots, n\} \;\;\mbox{and}\;\; K''=\{ i_1\},
\]
or
\[
K'= \{i_1\} \cup 
\{ n-m+1,\ldots, \widehat{l}, \ldots, n\} \;\;\mbox{and}\;\; K'' =\{ l\}
\]
where $n-m+1 \leq l \leq n$ and 
$l\notin \{ i_2, \ldots , i_m\}$.  
Let $S= \{n-m+1 , \ldots, n \}
\setminus \{  i_2, \ldots , i_m\}$.  By re-arranging the above equation,
we obtain 
\begin{align*}
\\[-3ex]
 [n-m+1, \dots, n]\qs{i_1\ldots i_m} = -\sum_{l\in S} \l
-q\r^\bullet\qs{i_1\; n-m+1\ldots \widehat{l} \ldots n}\qs{l \; i_2\ldots
i_m}. 
\end{align*}
Multiplying through by $ [n-m+1, \dots, n]^{-2}$ from the right, and using
Lemma~\ref{mincomm} gives
\begin{align*}
\\[-3ex]
\{i_1\ldots i_m\} = \sum_{l\in S} \pm\l
-q\r^{\bullet}\{i_1\; n-m+1\ldots \widehat{l} \ldots n\}
\{l \; i_2\ldots i_m\}.
\end{align*}
Now $\{ l, i_2, \ldots , i_m\} \cap \{ 1, \ldots , n-m\} = \{i_2, \dots,
i_{t-1}\}$ 
and so, by the inductive hypothesis, $\{l \; i_2\ldots i_m\} \in A$.
Clearly $\{i_1\; n-m+1\ldots \widehat{l} \ldots n\}\in A$, 
therefore $\{i_1\ldots i_m\}\in
A$. This completes the inductive step and the result follows. 
\hspace{2ex}\smallbox\end{proof}

\begin{theorem} \label{rhoiso}
There is an isomorphism 
\[
\rho:\oq(M_{m,n-m}) \goesto \dhom(\gqmn, [n-m+1, \dots,
n])
\]
 which is defined on generators by $\rho(X_{ij}) = 
\{j\; n-m+1 \ldots \widehat{n-i+1}\ldots n\}$, for $1\leq i \leq m$ and
$1 \leq j \leq n-m$. 
\end{theorem} 

\begin{proof} 
In order to 
show that $\rho$ is a homomorphism we have to show that the images of the
$X_{ij}$ under $\rho$ still obey the relevant commutation relations. We
will make repeated use of the anti-endomorphism $\Gamma\circ\tau$
defined before Lemma~\ref{gens}. 
There are four types of products to consider. 

(1)  Let $1\leq i<l\leq m$ and $1\leq j\leq n-m$.  Then
$
X_{ij}X_{lj} = q X_{lj}X_{ij},
$
and so we must show that
$
\rho\l X_{ij}\r\rho\l X_{lj}\r = q \rho\l X_{lj}\r \rho\l X_{ij}\r$. 
Let $ t = n +1-i$ and $ s = n+1-l$. Note that $s <t$, 
and consider the product
\[
  \qs{j \quad n-m+1 \ldots \widehat{t} \ldots n } \qs{j \quad n-m+1
\ldots \widehat{s}
\ldots n }
\]
in $\gqmn$.  We can think of this as a product in $\oq(M_{m+1})$ where
the rows are indexed by $1, \dots, m+1$ and the columns by $j, n-m+1,
\dots, n$. 
Apply the anti-endomorphism $\Gamma\circ\tau$ to the commutation relation
$X_{ m+1,s}X_{m+1, t } = qX_{m+1,t }X_{m+1,s }$ we obtain:
\begin{align*}
&\qs{j \;
n-m+1 \ldots \widehat{t} \ldots n } 
\qs{j \; n-m+1\ldots \widehat{s} \ldots n} \\[2ex]
&\hspace{22ex}= q\qs{j \; n-m+1 \ldots \widehat{s} \ldots n } \qs{j \;
n-m+1
\ldots \widehat{t}
\ldots n }.\\[-3ex]
\end{align*}

Multiplying through this equation on the right 
by $[n-m+1, \dots, n]^{-2}$ on each side and using Lemma~\ref{mincomm} gives
\begin{align*}
\\[-3ex]
&\{j \;
n-m+1 \ldots \widehat{t} \ldots n \} \{j \; n-m+1\ldots \widehat{s} \ldots
n\} \\[2ex]
&\hspace{22ex}= q\{j \; n-m+1 \ldots \widehat{s} \ldots n \} \{j \; n-m+1
\ldots \widehat{t}
\ldots n \};
\end{align*}
that is, $\rho\l X_{ij}\r\rho\l X_{lj}\r = q \rho\l X_{lj}\r \rho\l
X_{ij}\r$.   

(2)  Let $1\leq j<r\leq n-m$ and $1\leq i\leq m$.  Then 
$
X_{ij}X_{ir} = qX_{ir}X_{ij}.
$
Let $t=n+1-i$ and, as in (1), think of the product
\[
\qs{j \; n-m+1 \ldots \widehat{t} \ldots n}\qs{r\;  n-m+1 \ldots
\widehat{t} \ldots n}
\]
as sitting inside $\oq(M_{m+1})$ where the rows are indexed by $1,
\dots, m+1$ and the columns by $j, r, n-m+1, \dots \widehat{t}, \dots,
n$.   Then $\Gamma\circ\tau$ applied to the relation
$X_{m+1, j}X_{m+1, r } = qX_{m+1,r }X_{m+1,j } $ in $\oq(M_{m+1})$ 
gives us
\begin{align*}
&\qs{j \; n-m+1 \ldots \widehat{t} \ldots n }
\qs{r \; n-m+1\ldots \widehat{t} \ldots n } \\[2ex]
&\hspace{22ex}=q\qs{r \; n-m+1 \ldots \widehat{t} \ldots n } \qs{j \; n-m+1
\ldots \widehat{t}
\ldots n }.\\[-3ex]
\end{align*} 
Therefore, multiplying through this equation on the right 
by $[n-m+1, \dots, n]^{-2}$ and using Lemma~\ref{mincomm}, we get
\begin{align*}
\\[-3ex]
&\{j \; n-m+1 \ldots \widehat{t} \ldots n \}
\{r \; n-m+1\ldots \widehat{t} \ldots n \} \\[2ex]
&\hspace{22ex}=q\{r \; n-m+1 \ldots \widehat{t} \ldots n \} \{j \quad n-m+1
\ldots \widehat{t}
\ldots n \};
\end{align*}
that is, $\rho(X_{ij})\rho(X_{ir}) = q\rho(X_{ir})\rho(X_{ij})$

(3)  Let $1\leq i<l \leq m$, and $1 \leq j<r\leq n-m$.  Then
\[
X_{ij}X_{lr}=X_{lr}X_{ij} +\l q-q^{-1} \r X_{lj}X_{ir}.
\]
Let $t=n+1-i$ and $s=n+1-l$. Note that $n-m+1 \leq s < t \leq n$, and
that $j < r < s < t$. Consider the product
\[
\qs{ j\; n-m+1 \ldots \widehat{t} \ldots n}\qs{ r\;
n-m+1 \ldots \widehat{s} \ldots n}
\]
as a product in $\oq(M_{m+2})$, where the $m+2$ rows are indexed by $1,
\dots , m+2$ and the columns by $j, r, n-m+1, \dots , n$.  

The relation 
\[
\qs{13}\qs{24}  =  \qs{24}\qs{13} + \l q-q^{-1} \r \qs{14}\qs{23}
\]
that we calculated earlier for $\gqtwofour$ shows that, in
$\oq(M_{m+2})$,  
\[
[I\mid js][I\mid rt] = [I\mid rt][I\mid js] + (q - q^{-1}) [I\mid jt][I\mid rs]
\]
where $I = \{m+1,m+2\}$, since $j < r < s < t$. By applying 
the anti-endomorphism $\Gamma\circ\tau$ to this relation, we obtain 




\begin{align*}
&\qs{j\;n-m+1\ldots
\widehat{t}\ldots n}\qs{r\;n-m+1\ldots\widehat{s}\ldots
n} \\[2ex]
&\hspace{7ex}= \qs{r\;n-m+1\ldots\widehat{s}\ldots
n}\qs{j\;n-m+1\ldots
\widehat{t}\ldots n}\\[1ex]  
&\hspace{11ex}+\l q-q^{-1}\r\qs{j\;n-m+1\ldots
\widehat{s}\ldots n}\qs{r\;n-m+1\ldots
\widehat{t}\ldots n} 
\end{align*}
in $\gqmn$.
Multiplying through by $[n-m+1, \dots, n]^{-2}$ and using
Lemma~\ref{mincomm} 
we get
\begin{align*}
&\{j\;n-m+1\ldots
\widehat{t}\ldots n\}\{r\;n-m+1\ldots\widehat{s}\ldots
n\}\\[2ex]
&\hspace{7ex}= \{r\;n-m+1\ldots\widehat{s}\ldots
n\}\{j\;n-m+1\ldots
\widehat{t}\ldots n\}\\[1ex]
&\hspace{11ex}+\l q-q^{-1}\r\{j\;n-m+1\ldots
\widehat{s}\ldots n\}\{r\;n-m+1\ldots
\widehat{t}\ldots n\};
\end{align*}
that is, $\rho (X_{ij})\rho (X_{lr}) 
=\rho (X_{lr})\rho (X_{ij})+ (q-q^{-1})\rho
(X_{lj})\rho (X_{ir})$, as required. 
  
(4)  Let $1\leq i<l \leq m$ and $1 \leq j<r\leq n-m$.  Then
\[
 X_{ir}X_{lj}=X_{lj}X_{ir}.
\]
Let $t=n+1-i$ and  $s= n+1-l$ so that $n-m+1 \leq s <t \leq n$ and
$j<r<s<t$. Arguing as in (3), the relation $[23][14] = [14][23]$ in
$\gqtwofour$ produces, in $\oq(M_{m+2})$, the relation 
\[
[I\mid rs][I\mid jt] = 
[I\mid jt][I\mid rs].
\]
Applying $\Gamma\circ\tau$ to
this relation gives 
\begin{align*}
&\qs{ r\;n-m+1\ldots\widehat{t}\ldots
n}\qs{j\;n-m+1\ldots\widehat{s}\ldots n}\\[2ex]
&\hspace{22ex}= \qs{j\;n-m+1\ldots\widehat{s}\ldots
n}\qs{r\;n-m+1\ldots\widehat{t}\ldots n}.
\end{align*}
Multiplying through by $[n-m+1, \dots, n]^{-2}$ we get
\begin{align*}
&\{ r\;n-m+1\ldots\widehat{t}\ldots
n\}\{j\;n-m+1\ldots\widehat{s}\ldots n\}\\[2ex]
&\hspace{22ex}= \{j\;n-m+1\ldots\widehat{s}\ldots
n\}\{r\;n-m+1\ldots\widehat{t}\ldots n\};
\end{align*}
that is, $\rho(X_{ir})\rho(X_{lj}) =\rho(X_{lj})\rho(X_{ir})$, as
required. 

Thus, $\rho$ extends to  a homomorphism.  
The images of the generators under
$\rho$ generate $\dhom(\gqmn, [n-m+1, \dots, n])$, by
Lemma~\ref{gens}; so $\rho$ is surjective.  We show that $\rho$ is 
injective by comparing Gelfand-Kirillov dimensions.  If $\rho$ was not
injective, then $\gkdim(\dhom(\gqmn, [n-m+1, \dots, n]) <
\gkdim(\oq(M_{m,n-m})) = m(n-m)$, since $\oq(M_{m,n-m})$ is a domain.  
However,  by Corollary~\ref{dhomprops} and
Proposition~\ref{gragk}, we know that
$\gkdim(\dhom(\gqmn, [n-m+1, \dots, n]) = \gkdim(\gqmn) - 1 = 
m(n-m) +1 -1 = m(n-m)$.   Thus, $\rho$ is
injective and hence $\rho$ is an isomorphism. 
\hspace{2ex}\smallbox\end{proof}

\begin{corollary} \label{keyiso} 
Let $\phi$ be the automorphism of $\oq(M_{m,n-m})$ defined by $\phi(X_{ij}) =
q^{-1}X_{ij}$, for $1\leq i \leq m$ and $1\leq j\leq n-m$. Then 
\[
\oq(M_{m,n-m})[y,y^{-1};\phi] \goesto 
\gqmn\left[ [n-m+1, \dots, n]^{-1}\right]
\]
defined by $X_{ij} \mapsto \{j\, n-m+1 \dots \widehat{n+1-i} \dots
n\}$ and $y\mapsto [n-m+1, \dots, n]$ is an isomorphism of algebras.
\end{corollary}

\begin{proof} 
Recall from Lemma~\ref{sle} that there is an isomorphism
\[
\theta:\dhom(\gqmn, 
[n-m+1, \dots, n])[y,y^{-1};\sigma] \goesto \gqmn\left[
[n-m+1, \dots, n]^{-1}\right] 
\]
given by $y\mapsto [n-m+1, \dots, n]$ 
and $\{j\, n-m+1 \dots \widehat{t} \dots
n\} \mapsto \{j\, n-m+1 \dots \widehat{t} \dots     
n\}$, where $\sigma$ is the 
automorphism of $\dhom(\gqmn, [n-m+1, \dots, n])$ given
by conjugation by the quantum minor $[n-m+1, \dots, n]$. 
On the other hand, by Theorem~\ref{rhoiso}, there is an isomorphism 
$\rho:\oq(M_{m,n-m}) \goesto 
\dhom(\gqmn, [n-m+1, \dots, n])$, and it is
easy to see, by using Lemma~\ref{mincomm}, that the automorphism induced
in $\oq(M_{m,n-m})$ by $\sigma$ via 
$\rho$ is $\phi$. Thus, $\rho$ extends to an
isomorphism 
\[
\overline{\rho}: \oq(M_{m,n-m})[y,y^{-1};\phi] \goesto \dhom(\gqmn, 
[n-m+1, \dots, n])[y,y^{-1};\sigma] 
\]
such that $\overline{\rho}(y) =y$. Clearly, $\theta\circ\overline{\rho}$
is the desired isomorphism. 
\hspace{2ex}\smallbox\end{proof} \\

Note that in \cite{F2} Fioresi proves a restricted version of
Theorem~\ref{rhoiso}.  
More specifically, operating over the ring $K[q,q^{-1}]$, where
$K$ is algebraically closed of characteristic zero and $q$ is
transcendental over $K$, she shows that $\oqnn$ is isomorphic to the
subalgebra of $\gqntwon[\u^{-1}]$ generated by the elements $\{j \quad
n+1 \ldots \widehat{i} \ldots 2n\}$, but does not show that this
subalgebra is the dehomogenisation of $\gqntwon$ at $\u$. \\

%

\noindent{\bf Example}~~Let $S= \gqtwofour[{\qs{34}}^{-1}]$.  Then
$\dhom(\gqtwofour,\qs{34})= S_0$
and $S_0$ is generated by the elements
\[
\qs{12}\qs{34}^{-1},\quad \qs{13}\qs{34}^{-1},\quad
\qs{14}\qs{34}^{-1},\quad
\qs{23}\qs{34}^{-1},\quad \qs{24}\qs{34}^{-1}.
\]
Recall that $\{ ij\} = [ij]\qs{34}^{-1}$.  From the commutation
relations for $\gqtwofour$ given in the introduction, we 
can calculate the following commutation relations: 
\[
\{13\}\{23\} =q\{23\}\{13\};\quad
\{13\}\{14\}= q\{14\}\{13\};
\]
\[
\{13\}\{24\} =\{24\}\{13\} + \l q-q^{-1}\r
\{23\}\{14\};
\]
\[
\{14\}\{23\} =\{23\}\{14\}; \quad
\{14\}\{24\} = q\{24\}\{14\}; \quad \{23\}\{24\} = q\{24\}\{23\}
\]
and from the Quantum Pl\"{u}cker relation;
\[
\{12\} = \{13\}\{24\}-q\{23\}\{14\}.
\]
We can immediately see the correspondence (or we can use $\rho$ to
find the correspondence):
\[
\begin{array}{ccc}
\oq(M(2)) & \longleftrightarrow & S_0 \\
X_{11} & \longleftrightarrow & \{13\}\\
X_{12} & \longleftrightarrow & \{23\}\\
X_{21} & \longleftrightarrow & \{14\}\\
X_{22} & \longleftrightarrow & \{24\}\\
D_q & \longleftrightarrow & \{12\},
\end{array}
\]
and from Theorem~\ref{rhoiso}
\[
\dhom(\gqtwofour,[34]) \cong \oq(M(2)).
\]

\section{$\gqmn$ as coinvariants of $\oqslm$}

Recall that the 
{\em $m\times m$ quantum special linear group}, $\oqslm$, is defined
by $\oqslm:= \oqmm/\langle D_q -1\rangle$. 

	In this section we show that $\gqmn$ is the algebra of
coinvariants of a natural left coaction of $\oqslm$ on $\oqmn$.  There
is a natural epimorphism $\pi: \oqglm \goesto \oqslm$ which sends $D_q$
to $1$.  In order to distinguish generators in the various algebras, we
will often denote the canonical generators in $\oqnn$ by 
$X_{ij}$, in $\oqnm$ by $Y_{ij}$, in $\oqmn$ by $Z_{ij}$ and in $\oqglm$ 
by $T_{ij}$. Further, set $U_{ij}:= \pi (T_{ij}) \in \oqslm$.  Note
that both $\oqglm$ and $\oqslm$ are Hopf algebras. 

It is easy to check that one can define a morphism of algebras
satisfying the following rule:

\[
\lambda : \oqmn \goesto \oqglm \otimes \oqmn\quad Z_{ij} \mapsto
\sum_{k=1}^m\, T_{ik} \otimes Z_{kj}
\]
and that this induces a morphism of algebras  
\[
\Lambda : \oqmn \goesto \oqslm \otimes \oqmn, \quad Z_{ij} \mapsto
\sum_{k=1}^m\, U_{ik} \otimes Z_{kj}  
\]
where $\Lambda:= (\pi \otimes \id)\circ \lambda$. 

The morphisms $\lambda$ and $\Lambda$ endow $\oqmn$ with left comodule
algebra structures over $\oqglm$ and $\oqslm$, respectively.  Recall
that if $H$ is a Hopf algebra and $M$ is a left $H$-comodule via the
coaction $\gamma: M \goesto H\otimes M$ then $m \in M$ is a {\em
coinvariant} if $\gamma(m) = 1 \otimes m$.  In this section we show that
$\gqmn$ is the set of coinvariants of the $\oqslm$-comodule $\oqmn$
under the comodule map $\Lambda$.  In fact, this result is an easy
consequence of \cite[Theorem 6.6]{GLR}, once we have described the
set-up of that paper. 

The map $Y_{ij} \mapsto \sum_{k=1}^m\, Y_{ik} \otimes T_{kj}$ induces a
morphism of algebras $\rho: \oqnm \goesto \oqnm \otimes \oqglm$ which
endows ${\cal O}_q(M_{nm})$
with a right comodule algebra stucture over ${\cal
O}_q(GL_{m})$. Let
$\oqv$ denote the algebra $\oqnm \otimes \oqmn$. The coactions $\lambda$
and $\rho$ defined above can be combined to give a left comodule
structure on $\oqv$ which we denote by $\gamma$. To be precise, 
\[
\gamma : \oqv \goesto \oqglm \otimes \oqv
\]
is given by the rule 
\[
\gamma (a \otimes b):= \sum_{(a),(b)}  S(a_{1})b_{-1} \otimes a_{0}
\otimes b_{0}
\]
for $a\in \oqnm$ and $b \in \oqmn$, 
where $\lambda(b) = \sum_{(b)}\, b_{-1} \ot b_0$ and $\rho(a) =
\sum_{(a)}\, a_0\ot a_1$. Here, 
we are using the Sweedler
notation and $S$ is the antipode of $\oqglm$. In turn, this coaction
induces a coaction $\Gamma :\oqv \goesto \oqslm \otimes \oqv$ given by
$\Gamma:= (\pi \otimes \id)\circ \gamma$; so that 
\[
\Gamma (a \otimes b):= \sum_{(a),(b)}  \pi(S(a_{1})b_{-1}) 
\otimes a_{0} \otimes b_{0}.
\]
The main results of \cite{GLR}
identify the coinvariants of the coactions $\gamma$ and $\Gamma$.  
In particular, 
Theorem 6.6 of \cite{GLR} identifies the coinvariants of the coaction
$\Gamma$ in the following way. 
There is a morphism of algebras $\mu:\oqnn \goesto \oqv =
\oqnm \otimes \oqmn$ given by $X_{ij} \mapsto \sum_{k =1}^m\, Y_{ik}
\otimes Z_{kj}$. Let $R$
denote $\mu(\oqnn)$. It is proved in \cite{GL} that $R\cong
\oqnn/I$, where $I$ is the ideal generated by the $(m+1)\times (m+1)$
quantum minors of $\oqnn$. We have the following theorem.

\begin{theorem}\cite[Theorem 6.6]{GLR} \label{glr6.6}
Let $G_1$ and $G_2$ denote the respective grassmannian subalgebras 
of $\oqnm$ and $\oqmn$ generated by all the     
$m\times m$ quantum minors. The set of $\Gamma$-coinvariants in
$\oqv= \oqnm\otimes\oqmn$ is the subalgebra generated by
$G_1\otimes G_2$ and $R$. More precisely,
$$\left( \oqnm\otimes\oqmn \right)^{\co\, \oqslm}= \left(
G_1\otimes G_2 \right) \cdot R.$$ 
\end{theorem} 

The result we are aiming for follows easily from this.

\begin{theorem}
\[
\left(\oqmn\right)^{\co\,\oqslm} = \gqmn.
\]
\end{theorem} 

\begin{proof} 
It is easily seen that there is  a
commutative diagram  
\[
\begin{array}{ccc}
\oqmn &
\stackrel{i}{\longrightarrow} &   
\oqnm \otimes \oqmn \cr
&& \cr
\Lambda \;\downarrow \qquad
& & \Gamma\;\downarrow \qquad \cr
&&\cr {\cal O}_q(SL_m)\otimes\oqmn &
\stackrel{{\rm id} \otimes i}{\longrightarrow} &
{\cal O}_q(SL_m)\otimes\oqnm \otimes \oqmn
\cr
\end{array}
\]
where $i$ is the canonical injection.  Moreover, let $j : 
\oqnm\otimes \oqmn \longrightarrow \oqmn$ 
be the canonical projection; that is, 
\[
j : \oqnm \otimes \oqmn
\stackrel{p\otimes {\rm id}}{\longrightarrow}
k \otimes \oqmn \cong \oqmn 
\]
where $p$ is the projection modulo the irrelevant ideal of ${\cal
O}_q(M_{nm})$.  Clearly, we have that $j\circ i ={\rm id}$.  We see from
the above commutative diagram that, if $b \in {\cal O}_q(M_{mn})$ is a
$\Lambda$-coinvariant, then $i(b)=1 \otimes b$ is a $\Gamma
$-coinvariant. 
Thus, it follows from Theorem~\ref{glr6.6} that $1
\otimes b \in (G_1 \otimes G_2).R $.  
Hence, 
$b=j(1 \otimes b) \in 
j(G_1 \otimes G_2)j(R)$.  
Clearly, $j(R) \subseteq k$ and
$j((G_1 \otimes G_2)) \subseteq G_2$; and so $b
\in G_2=\gqmn$.  This shows that $  {\cal O}_q(M_{mn})^{{\rm co} {\cal
O}_q(SL_m)}\sse \gqmn$.  
Since it is
clear that an $m\times m$ quantum 
minor of ${\cal O}_q(M_{mn})$ is a $\Lambda$-coinvariant, 
the converse inclusion follows from the fact that
$\Lambda$ is a morphism of algebras. 
\hspace{2ex}\smallbox\end{proof}\\

Note that Fioresi and Hacon, \cite{FH}, have a version of this result,
with the usual restrictions as described earlier in this paper.

\section{$\gqmn$ is a maximal order}

Let $R$ be a noetherian domain with division ring of fractions $Q$.
Then $R$ is said to be a {\bf maximal order} in $Q$ if the following   
condition is satisfied: if $T$ is a ring such that $R\subseteq T 
\subseteq Q$ and such that there exist nonzero elements $a, b \in R$
with $aTb\sse R$, then $T=R$.  This condition is the natural
noncommutative analogue of normality for commutative domains, see, for
example, \cite[Section 5.1]{McCR}. 

Recall that an element $d$ in a ring $R$ is said to be {\em left
regular} if $rd = 0$ implies that $r=0$ for $r\in R$.  The following is
a general result that we will be able to apply to show that the quantum
Grassmannian $\gqmn$ is a maximal order. 

\begin{proposition} 
Suppose that $R$ is a noetherian domain with
division ring of fractions $Q$. Suppose that $a, b \in R$ are nonzero 
normal elements such that $R[a^{-1}]$ and $R[b^{-1}]$ are both maximal
orders, that  $b$ is 
left regular modulo $aR$ and that $ab = \lambda ba$ for some central
unit $\lambda \in R$. 
Then $R$ is a maximal order. 
\label{propmax}
\end{proposition} 

\begin{proof} 
First, we show that 
$R[a^{-1}] \cap R[b^{-1}] = R$. 
Suppose that this is not the
case, and choose $q \in R[a^{-1}] \cap R[b^{-1}]\; \backslash R$. Write
$q = ra^{-d} = sb^{-e}$ with $d, e\geq 1$ and $r \in R\backslash Ra$, 
$s\in R\backslash Rb$.  Cross
multiply to get $rb^e = \lambda^{\bullet} sa^d$ 
(remember that $ab = \lambda ba$).
Since $b$ is left regular modulo $aR$, 
this gives $r\in Ra$, a contradiction. Thus,  $R[a^{-1}] \cap R[b^{-1}]
= R$.

Now, to show that $R$ is a maximal order, it is enough to show that if
$J$ is a nonzero ideal of $R$ and $q\in Q$ with either $qJ \subseteq J$ 
or $Jq \subseteq J$ then $q \in R$, \cite[Proposition 5.1.4]{McCR}. 
Suppose, without loss of generality, that $qJ \subseteq J$. 
By assumption, $S :=
R[a^{-1}]$ and $T:= R[b^{-1}]$ are maximal orders. 
Also, $SJ = JS$ is an ideal of $S$ and $TJ = JT$ is an
ideal of $T$.  We have $qJS \subseteq JS$ and so $q \in S$.  Similarly,
$q\in T$.  Thus, $q \in S\cap T = R$; and so $R$ is a maximal order. 
\hspace{2ex}\smallbox\end{proof} 

\begin{theorem}$\gqmn$ is a maximal order. 
\end{theorem} 

\begin{proof} We will apply the previous result to 
$R:= \gqmn$ with $a:= [1, \dots, m]$ and $b:= [n-m
+1, \dots, n]$. 
Observe that $b$ is normal by Lemma~\ref{mincomm} 
and that $a$ is normal by
Corollary~\ref{commr}.
Note that $ab = \qdot ba$, by Lemma~\ref{mincomm}. 
First we observe that $b$ is left regular 
modulo $aR$. The reason is that since $a$ is the 
minimal minor in the preferred ordering, a basis for $aR$ is given by
preferred products that start with $a$. If $r \in R$ is such that $rb
\in aR$, then when we write $r$ as a linear combination of preferred
products then multiplying each preferred product that occurs by $b$ on
the right 
still gives a preferred product, since $b$ is the maximal element with
respect to the preferred order. Thus, since $rb \in aR$ each of these
preferred products must begin with $a$, and so the original ones also
begin with $a$, hence $r\in aR$. 

In Corollary~\ref{keyiso}, 
we have shown that 
$R[b^{-1}] \cong {\cal O}_q(M_{m,n-m})[y,y^{-1}; \phi]$ and so
$R[b^{-1}]$ is a maximal order (\cite[V. Proposition 2.5, IV.
Proposition 2.1]{MR}).  Also $R[a^{-1}]$ is 
a maximal order by using the isomorphism $\delta$ introduced in
Section 1 and the fact that $R[b^{-1}]$ is a maximal order.  

Thus, the hypotheses of Proposition~\ref{propmax} are satisfied, and we
deduce that $\gqmn$ is a maximal order. 
\hspace{2ex}\smallbox\end{proof}




\noindent A. C. Kelly, T. H. Lenagan:
\\School of Mathematics, University of Edinburgh,
\\ James Clerk Maxwell Building, King's Buildings, Mayfield Road,
\\Edinburgh EH9 3JZ, Scotland
\\E-mail: ACK@bosinternet.com, tom@maths.ed.ac.uk
\\
\\
L. Rigal: 
\\Universit\'e Jean Monnet
(Saint-\'Etienne), Facult\'e des Sciences et
\\Techniques, D\'e\-par\-te\-ment de Math\'ematiques
\\ 23 rue du Docteur Paul Michelon
\\42023 Saint-\'Etienne C\'edex 2,\\France
\\E-mail: Laurent.Rigal@univ-st-etienne.fr



\begin{thebibliography}{99}

\bibitem{ATV} M Artin, J Tate and M van den Bergh, {\it Some algebras
associated to automorphisms of elliptic curves}, The Grothendieck
Festschrift, Vol.  I, 33-85, Progr.  Math., {\bf 86}, Birkh\"auser
Boston, Boston, MA, 1990.


\bibitem{BV} W Bruns and U Vetter, {\it Determinantal rings}, 
Springer Lecture Notes in
Mathematics, 1327, Springer-Verlag, Berlin,
1988.

\bibitem{F} R Fioresi,  {\it Quantum deformation of the Grassmannian
manifold},  J. Algebra {\bf 214} (1999), 418-447.

\bibitem{F2} R Fioresi, {\it A deformation of the big cell inside the
Grassmannian manifold $G(r,n)$}, Rev. Math. Phys. {\bf 11} 
(1999), 25-40.

\bibitem{FH} R Fioresi and C Hacon, {\it Quantum coinvariant theory for
the quantum special linear group and quantum Schubert varieties}, 
J. Algebra {\bf 242} (2001), 433-446.

\bibitem{GL} K R Goodearl and T H Lenagan. {\it Quantum determinantal  
ideals}, Duke Mathematical Journal  {\bf 103} 165-190,  2000.

\bibitem{GL2} K R Goodearl and T H Lenagan. {\it Winding-invariant prime
ideals in quantum $3\times 3$ matrices}, to appear in Journal of Algebra, 
preprint available at math.QA/0112051. 

\bibitem{GLR} K R Goodearl, T H Lenagan and L Rigal, {\it The first
fundamental theorem of coinvariant theory for the quantum general linear
group}, Publ. RIMS (Kyoto) {\bf 36} (2000), 269-296. 

\bibitem{KL}G R Krause and T H Lenagan. {\it Growth of algebras and
Gelfand-Kirillov dimension}, Revised edition.
Graduate Studies in Mathematics, 22. American Mathematical Society,
Providence, RI, 2000.

\bibitem{MR} G Maury and J Raynaud, {\it Ordres maximaux au sens de K
Asano}, Springer Lecture Notes in Mathematics Vol 808, Springer-Verlag,
Berlin, 1980.


\bibitem{McCR} J C McConnell and J C Robson. {\it Noncommutative   
Noetherian Rings.} Wiley, Chichester, 1987.

\bibitem{NYM} M Noumi, H Yamada and K Mimachi, {\it
Finite-dimensional representations of the quantum group ${\rm
GL}_q(n;\mc)$ and the zonal spherical functions on $U_q(n-1)\backslash
U_q(n)$}, Japanese J. Math {\bf 19} (1993), 31-80. 

\bibitem{PW}  B Parshall and J Wang. {\it Quantum linear groups.} Mem.
Amer. Math. Soc {\bf 89} (1991), no. 439.

\bibitem{Z} J J Zhang, {\it Connected graded Gorenstein algebras with
enough normal elements}, J. Algebra {\bf 189} (1997), 390-405.

\end{thebibliography}
\end{document}